\documentclass[12pt]{article}
\usepackage{amsmath,epsfig,amssymb,float,graphicx,amsfonts}
\usepackage{xspace,colortbl}

\definecolor{color1}{rgb}{0.3,0.5,1}

\title{On the Stochastic Rank of Metric Functions}

\author{Nikolay Balov}

\begin{document}

\maketitle

\begin{abstract}
For a class of integral operators with kernels metric functions on manifold we find some necessary and sufficient conditions to have finite rank. 
The problem we pose has a stochastic nature and boils down to the following alternative question. For a random sample of discrete points, what will be the probability the symmetric matrix of pairwise distances to have full rank? When the metric is an analytic function, the question finds full and satisfactory answer. 
As an important application, we consider a class of tensor systems of equations 
formulating the problem of recovering a manifold distribution from its covariance field and solve this problem for representing manifolds such as Euclidean space and unit sphere. 
\end{abstract}

\section{Problem formulation and motivation}

We start with the classical Fredholm integral equation of the first kind 
\begin{equation}\label{eq:rankmetric_fredholm}
\int_V \psi(x,y) f(y)dy = g(x), \textrm{ }x\in U,
\end{equation}
where $U$ and $V$ are open sets in $\mathbb{R}^n$ and $f:U\to\mathbb{R}$, $g:V\to\mathbb{R}$, 
and $\psi:U\times V\to\mathbb{R}$ are some functions. Depending on the domains more 
conditions on $f$, $g$ and $\psi$ may be necessary for the correct formulation of (\ref{eq:rankmetric_fredholm}). 
Let $\{x_i\}_{i=1}^k$ and $\{y_i\}_{i=1}^k$ be two discrete samples of points chosen 
by uniform distributions on $U$ and $V$ respectively. 
Then equation (\ref{eq:rankmetric_fredholm}) can be discretized by the matrix-quadrature method 
\begin{equation}\label{eq:rankmetric_fredholm_discrete}
\sum_{j=1}^k \psi(x_i, y_j) f(y_j) = g(x_i), \textrm{ }i=1,...,k.
\end{equation}
In general, the inverse problem of solving (\ref{eq:rankmetric_fredholm}) for $f$, 
is often ill-posed and approximation based on (\ref{eq:rankmetric_fredholm_discrete}) 
will eventually result in increasingly unstable solution as $k$ increases. 
Here, however, we are interested in the first potential obstacle to solve (\ref{eq:rankmetric_fredholm_discrete}) - 
the matrix $\Psi := \{\psi(x_i,y_j)\}_{i,j=1,1}^{k,k}$ may not be of full rank. 
The problem is stochastic one  
for the points $x_i$'s and $y_j$'s are chosen in random fashion.
In section 2 we investigated it and find some conditions for the kernel $\psi$, 
which guarantee that for any $k$, $\Psi$ has full rank with probability one.
We also show some necessary and sufficient conditions 
for analytic kernels $\psi$ to be of finite rank.

A further generalization of equation (\ref{eq:rankmetric_fredholm}) takes $f$ to be a function 
on n-manifold M and $g$ and $\psi$ to be linear operator fields on M. 
For a point $p\in M$ with $M_p$ we denote the tangent space at $p$ and 
with $T_1^1(M_p)$, the vector space of (1,1)-tensors (linear operators) on $M_p$. 
Let $\mu$ be a measure on M as for example the volume measure $V(p)$ on Riemannian manifolds. 
Consider the equation 
\begin{equation}\label{eq:rankmetric_inverse_tensor_field}
\int_M Y(p,q) f(q) d\mu(q) = C(p), \textrm{ }p\in M,
\end{equation}
such that $Y(p,.)\in T_1^1(M_p)$ and $C(p)\in T_1^1(M_p)$. 
The inverse problem here is finding $f$ for given fields $Y$ and $C$.
If we know that (\ref{eq:rankmetric_inverse_tensor_field}) has a unique solution for $f$ then it 
can be found by solving  
\begin{equation}\label{eq:rankmetric_inverse_trace_field}
\int_M tr(Y(p,q)) f(q) d\mu(q) = tr(C(p)),
\end{equation}
an equation of type (\ref{eq:rankmetric_fredholm}). 

The importance of the class (\ref{eq:rankmetric_inverse_tensor_field}) of tensor equations is 
that it contains the problem of recovering a distribution from its covariance field. 
Next, we briefly pose this problem, while more details one can find in \cite{balov-covf}.

Let M be a Riemannian manifold with metric tensor $G$. For any $p\in M$, $G(p)\in T^2(M_p)$ is a co-variant 2-tensor.
Let $Exp_p: M_p \to M$ be the exponential map at $p$ and 
$\mathcal{U}(p)\subset M$ be the maximal normal neighborhood of $p$, 
where $Exp_p^{-1}$ is well defined. 
Note that since $Exp_p^{-1}q \in M_p$, $(Exp_p^{-1}q)(Exp_p^{-1}q)^T \in T_2(M_p)$, a contra-variant 2-tensor.
For a density function  $f\ge 0$ on M, the covariance operator field of $f$ is 
$G\Sigma:M \to T_1^1(M)$, such that for any $p\in M$
\begin{equation}\label{eq:rankmetric_covfield}
G\Sigma(p) := \int_{\mathcal{U}(p)} G(p)(Exp_p^{-1}q)(Exp_p^{-1}q)^T f(p)dV(p).
\end{equation}
The problem of distribution recovering is of type (\ref{eq:rankmetric_inverse_tensor_field})
if we take $\mu=V$, $C=G\Sigma$ and $Y(p,q) = G(p)(Exp_p^{-1}q)(Exp_p^{-1}q)^T$.
Note that 
$$tr(Y(p,q)) = tr(G(p)(Exp_p^{-1}q)(Exp_p^{-1}q)^T) = d^2(p,q),$$ is the square geodesic distance on M.
Thus equation (\ref{eq:rankmetric_inverse_trace_field}), specifically, is 
\begin{equation}\label{eq:rankmetric_sqdist_kernel}
\int_{\mathcal{U}(p)} d^2(p,q) f(q) dV(q) = g(p).
\end{equation}

In the context of problem (\ref{eq:rankmetric_sqdist_kernel}) we are interested in finding the rank of an integral operator of the form 
$L_{\psi}: f \mapsto \int \psi(p,q)f(q)dV(q),$ 
where $\psi(p,q) = d^2(p,q)$ is a square distance function on M. In particular, 
in section 3 we study 
the rank of the Euclidean metric $d(p,q) = ||p-q||$ and 
the rank of the standard metric on the unit sphere $\mathbb{S}^n$, $d(p,q) = \cos^{-1}(<p, q>)$, 
and show that while the Euclidean metric is of finite rank, the metric on the sphere is not. 
The last fact can be re-phrased as follows. For any discrete sample of points on a sphere, 
the square matrix of pairwise distances is non-singular with probability one.

The problem of establishing the non-singularity of the kernel of operator (\ref{eq:rankmetric_fredholm})
is important in the context of more general statistical inverse problems on manifolds, 
as considered in \cite{cavalier-tsybakov}, \cite{kim-invmani}, \cite{koo-optimalrates} and \cite{ruymgaart-inverse}.
Let $g(p) = f(p) + \epsilon$, where $\epsilon$ is a mean zero random variable with small variance, 
be a model with unknown regression function $f$.
For a kernel $\psi$, the inverse problem with random noise is formulated as estimation of $L_{\psi}f$ 
from observations $(p_i,g_i)$. 
In \cite{cavalier-tsybakov}, Cavalier and Tsybakov estimate $f$ from the model 
$g = L_{\psi} f + \epsilon$, which is a noised version of (\ref{eq:rankmetric_fredholm}).
Usually points $p_i$ are assumed uniformly distributed on M. 
If the kernel $\psi$ has a finite rank,  
then these problems are necessarily ill-posed, 
since the operator $L_{\psi}$ is not invertible. 

Statistical inverse estimation is addressed in details in 
\cite{koo-optimalrates}, \cite{ruymgaart-inverse}, \cite{cavalier-tsybakov}, \cite{ruymgaart-short} and \cite{kim-invmani}.
In addition, in \cite{cavalier-inverse} different alternatives for regularization of ill-posed linear problems are  discussed. 
In our final remarks we also briefly elaborate on this important aspect of all inverse problems. 
However, instead of regularizing a fixed kernel, we alleviate the ill-conditioning by choosing a kernel with 
better conditional number from a class of similar kernels. This is a completely different kind of approach. 
At the end, we show some simulation results supporting our proposal.

It is a common assumption for $L_{\psi}$ to be Hermitian and even compact operator, although 
Mair and Ruymgaart, \cite{mair-ruymgaart}, and Cavalier, \cite{cavalier-noncompact}, relax the assumptions to Hermitian and positive definiteness. 
There are symmetric operators however, as the covariance one shown above, for which even these assumptions are too strong 
and this makes relevant the subject of this study - 
the rank of a general operator with analytic kernel.


\newtheorem{rankmetric_remark1}{Remark}
\begin{rankmetric_remark1}
Provided that there is a local chart $(\tilde U, \phi)$, $\tilde U\subset \mathbb{R}^n$ on M 
such that $\phi(\tilde U)$ contains the support of $f$, $Y$ and $C$ , equations 
(\ref{eq:rankmetric_inverse_tensor_field}), (\ref{eq:rankmetric_covfield}) and (\ref{eq:rankmetric_inverse_trace_field}) 
can be formulated as integral equations in Euclidean space $\mathbb{R}^n$ by 
substituting $p$ and $q$ with $\phi(x)$ and $\phi(y)$ respectively.
In fact, the above requirement is not a strong one for complete Riemannian manifolds where $V(M\backslash \mathcal{U}(p)) = 0$.
That is why, for the sake of simplicity, to the end of this exposition we will work in Euclidean settings.
\end{rankmetric_remark1}

\section{Rank of bi-variate analytic functions}

Our goal here is to find some necessary and sufficient conditions for an integral kernel to be of finite rank. 
We restrict to real analytical kernels, i.e. bi-variate real analytical functions.
The exposition uses elementary functional analysis techniques with some textbook facts presented for the sake of consistency.

Let $U$ be an open subset of $\mathbb{R}^n$. 
A collection of functions $f_s:U \to \mathbb{R}$ is said to be linear independent in $U$ if 
$\sum_s \alpha_s f_s(x) = 0$, for almost all (by Lebesgue measure) $x\in U$ only if $\alpha_s = 0$ for all $s$.
\newtheorem{rankmetric_lemma1}{Lemma}
\begin{rankmetric_lemma1}
If functions $f_1(x)$, ...,$f_k(x)$ are linearly independent in $U$, then there exist 
$x_1,...,x_k$, such that $rank(\{f_i(x_j)\}_{i=1,j=1}^{k,k}) = k$.
\label{lemma:linindep_functions}
\end{rankmetric_lemma1}
\begin{proof} Since $f_1(x)$ is not identically zero, there is $x_1\in U$ such that $f_1(x_1) \ne 0$. 
Determinant 
$$
det \left( \begin{array}{cc}
f_1(x_1) & f_2(x_1) \\
f_1(x) & f_2(x)
\end{array} \right) = f_1(x_1)f_2(x) - f_2(x_1)f_1(x)
$$
can not vanish for all $x\in U$ or $f_1$ and $f_2$ would be linear dependent. Therefore, there is $x_2$ 
such that $det(\{f_i(x_j)\}_{i=1,j=1}^{2,2}) \ne 0$. 
The selection process can be extended to find $x_1$,...,$x_k$ in $U$, such that $det(\{f_i(x_j)\}_{i=1,j=1}^{k,k}) \ne 0$. 
At the last step, defining $A_1:=det(\{f_i(x_j)\}_{i=1,j=1}^{k-1,k-1}) \ne 0$, 
we choose $x_k\in U$  such that for $x=x_k$ 
$$
det \left( \begin{array}{cccc}
f_1(x_1) & f_2(x_1) & ... & f_{k}(x_1) \\
... & ... & ... & ... \\
f_1(x_{k-1}) & f_2(x_{k-1}) & ... & f_{k}(x_{k-1}) \\
f_1(x) & f_2(x) & ... & f_{k}(x) \\
\end{array} \right) = \sum_{i=1}^{k} A_i f_i(x) \ne 0.
$$
\qquad\end{proof} 
\newtheorem{rankmetric_def1}{Definition}
\begin{rankmetric_def1}
We say that function $\psi:U\times V \to \mathbb{R}$, $U,V\subset\mathbb{R}^n$ 
has a rank $k$ and write $rank(\psi) = k$ if 
for any $m\in\mathbb{N}$, $x_i\in U$ and $y_i \in V$, i,j=1,...,m, 
$$
rank( \{\psi(x_i, y_i)\}_{i=1,j=1}^{m,m} ) \le k, 
$$
and $k$ is the smallest number with this property. 
\label{rankmetric:rank_function}
\end{rankmetric_def1}

Let $\psi:U\times V \to \mathbb{R}$, 
be a smooth bi-variate function, a fact we denote with $\psi\in C^{\infty}(U\times V)$.
We say that $\psi$ is analytic in $V$ about a point $p\in V$ if 
$\psi$ has Tailor expansion in $y\in V$   
$$
\psi(x,y) = \sum_{m=0}^{\infty} \sum_{s:\sum_{i=1}^n s_i = m} c_{s_1...s_n}(x) (y_1-p_1)^{s_1}...(y_n-p_n)^{s_n}, 
$$
which for the sake of brevity we will write as 
$$
\psi(x,y) = \sum_{[s]=0}^{\infty} c_s(x) (y-p)^s, 
$$
where $s=(s_1,...,s_n)$ is a multi-index and $[s]=s_1+...+s_n$. 

We say that the space spanned by the functions $c_s$ have a finite basis if there exist a number $m$ and 
functions $f_1$, ..., $f_m$, such that $c_s$'s are linear combinations of $f_l$'s, 
i.e. $c_s\in span\{f_1,...,f_m\}$ for all $s$.

Next result gives some necessary and 
sufficient conditions for $\psi$ to have a finite rank.
\newtheorem{rankmetric_th1}{Proposition}
\begin{rankmetric_th1}
For a function $\psi\in C^{\infty}(U\times V)$ that is analytic in $V$ about a point $p\in V$ 
the following three conditions are equivalent
\begin{enumerate}
\item[(1)] $rank(\psi) = k$.
\item[(2)] $\psi(x,y) = \sum_{j=1}^k f_j(x)g_j(y)$ for linearly independent 
functions $f_j\in C^{\infty}(U)$ and $g_j\in C^{\infty}(V)$, j=1,...,k.
\item[(3)] There are linearly independent functions $f_j\in C^{\infty}(U)$,j=1,...,k, such that 
		all Taylor functions $c_s$ are their linear combinations in $U$, i.e. 
		$c_s \in span\{f_1,...,f_k\}$.
\end{enumerate}
\label{theorem:rank_equivals}
\end{rankmetric_th1}
\begin{proof} Without loss of generality we assume that $p$ is the origin $0$. 
We show that $(3) \Rightarrow (2) \Rightarrow (1) \Rightarrow (3)$.

Let condition $(3)$ hold. Then for any $s$, $c_s(x) = \sum_{i=1}^k \beta_s^i f_i(x)$ and 
$\psi(x,y) = \sum_{i=1}^k (\sum_s \beta_s^i y^s) f_i(x)$, provided that 
for any $i$ and $y$, $\sum_s \beta_s^i y^s$ converges. 
Suppose the contrary, that there exists $i_0$ and $y_0\in V $ such that $\sum_s \beta_s^{i_0} y_0^s$ does not converge.
Choose $x_j\in U$,j=1,...,k, such that $rank(A := \{f_i(x_j)\}_{j=1,i=1}^{k,k}) = k$.
Let $||A^{-1}|| > 0$ denote the norm of matrix $A^{-1}$. Fix a number $\epsilon > 0$. 
Since all $\sum_s c_s(x_j) y_0^s$ converge, there is $N$ such that 
$$
|\sum_{[s]=N}^M (\sum_{i=1}^k \beta_s^i f_i(x_j)) y_0^s| < \epsilon / ||A^{-1}||,
$$ 
for any $j$ and $M\ge N$. Also by the assumption, there is $M > N$, such that 
$$
|\sum_{[s]=N}^M \beta_s^{i_0} y_0^s| \ge \epsilon.
$$ 
Define k-vectors  
$
z=(z_1,...,z_k)^T 
\textrm{ and } 
w=(w_1, ..., w_k)^T,
$ 
where 
$$
z_i := \sum_{[s]=N}^M \beta_s^{i} y_0^s
\textrm{ and } 
w_j := \sum_{[s]=N}^M (\sum_{i=1}^k \beta_s^{i} f_i(x_j))y_0^s.
$$ 
The system $A z = w$ can be solved for $z$, $z = A^{-1}w$. 
Then $||z|| \le ||A^{-1}||\textrm{ }||w|| < \epsilon$, 
which contradicts $||z|| \ge |z_{i_0}| \ge \epsilon$.
Therefore, the initial assumption is false and 
$$
g_i(y) := \sum_s \beta_s^i y^s,
$$
are well defined functions in $U$. Thus condition (2) holds.

That (1) follows from (2) is obvious.
Indeed, for $x_i\in U, y_i\in V$, i=1,...,k, define 
k-vectors $v_j = (g_j(y_1),...,g_j(y_k))^T$ and 
 $\psi_i = (\psi(x_i,y_1),...,\psi(x_i,y_k))^T$. 
Then for the column-vectors of 
matrix $\{\Psi=\psi(x_i,y_j)\}_{i,j=1}^{k,k}$ we have $\psi_i = \sum_{j=1}^k f_j(x_i) v_j$ which shows 
$rank(\Psi)=k$.

Next we show $(1) \Rightarrow (3)$. Let $rank(\psi) = k$. For fixed $x_i \in U$, i=1,...,k+1, consider the set 
$$\psi(x_1, y),...,\psi(x_{k+1}, y)$$ of functions of $y$. 
Since no full rank matrix $\{\psi(x_i,y_j)\}_{i=1,j=1}^{k+1,k+1}$ exists, by lemma \ref{lemma:linindep_functions}, 
$\{\psi(x_i, y)\}_{i=1}^{k+1}$ are not linear independent.
Thus there are functions $\alpha_i(x), x=(x_1,...,x_{k+1})$, not all zero, such that 
$$
\sum_{i=1}^{k+1} \alpha_i(x) \psi(x_i, y) = 0, \forall y\in V.
$$
By expanding $\psi$ in $y$ we obtain
$$
\sum_{s} \sum_{i=1}^{k+1} \alpha_i(x)c_s(x_i) y^s = 0, \forall y\in V,
$$
equivalent to 
$$
\sum_{i=1}^{k+1} \alpha_i(x) c_s(x_i) = 0, \forall s,
$$
and this is true for any $x_i\in U$, i=1,...,k+1.
Again by lemma \ref{lemma:linindep_functions}, no set of $k+1$ functions $c_s$ is linear independent.
Therefore, there exists a set of functions $f_1(x)$,...,$f_k(x)$, such that for any $s$, 
$c_s(x) = \sum_{i=1}^k \beta_s^i f_i(x)$, i.e. $c_s\in span\{f_1,...,f_k\}$.
Also, no less than $k$ such functions exists, because otherwise, following 
$(3)\Rightarrow (2) \Rightarrow (1)$, we would have that $rank(\psi) < k$.
Thus condition (3) holds.
\qquad\end{proof}

Note that if $\psi$ is symmetric and representation (2) holds, then $g_i$ are linear combinations of $f_i$.
Indeed, if we choose $x_j$,j=1,...,k, such that 
$
rank(\{f_i(x_j)\}_{i=1,j=1}^{k,k}) = k,
$ 
then we can solve the system 
$$
\sum_{i=1}^k f_i(x_j)g_i(y) = \sum_{i=1}^k f_i(y)g_i(x_j), j=1,...,k
$$
to obtain 
$
g_i(y) = \sum_{i=1}^k \gamma_{ij} f_j(y).
$

An immediate corollary of Theorem \ref{theorem:rank_equivals} is that for $\psi$ as assumed there, 
which has a finite rank, $rank(\psi) = k$, the integral operator 
$$
L_{\psi}: f \mapsto \int_V \psi(x,y)f(y)dy, \textrm{ } f:U\to \mathbb{R}, 
$$
is of finite rank.

\newtheorem{rankmetric_ex1}{Example}
\begin{rankmetric_ex1}\label{ex:rank_Rn}
Consider the Euclidean metric in $\mathbb{R}^n$, 
$
\psi(x,y) = \sum_{s=1}^n (x^s-y^s)^2,
$
where $x^s$ are the components of $x\in \mathbb{R}^n$.
Because of the global representation 
$$
\psi(x,y) = \sum_{s=1}^n (x^s)^2 * 1 - \underbrace{2x^1 y^1 ... - 2 x^n y^n}_{n} + 1 * \sum_{s=1}^n (y^s)^2,
$$
by Theorem \ref{theorem:rank_equivals}, $rank(\psi) = n+2$, i.e. 
the Euclidean metric has a finite rank of 2 more than the number of dimensions.
\end{rankmetric_ex1}

\newtheorem{rankmetric_def2}[rankmetric_def1]{Definition}
\begin{rankmetric_def2}
A collection of functions $f_s:U \to \mathbb{R}$ is said to be locally linear independent in $U$ if 
for any open subset $V$ of $U$ we have that 
$\sum_s \alpha_s f_s(x) = 0$, for almost all $x\in V$ only if $\alpha_s = 0$ for all $s$.
\label{def:local_linear_independence}
\end{rankmetric_def2}
Note that local linear independence is stronger that linear independence and the former implies the latter. 
\newtheorem{rankmetric_lemma2}[rankmetric_lemma1]{Lemma}
\begin{rankmetric_lemma2}
If functions $f_1(x)$, ...,$f_k(x)$ are locally linearly independent in $U$ and 
$U_i\subset U$,i=1,...,k are open subsets, then there exist 
$x_1,...,x_k$, such that $x_i\in U_i$ and \\ {$rank(\{f_i(x_j)\}_{i=1,j=1}^{k,k}) = k$}.
\label{lemma:locally_linindep_functions}
\end{rankmetric_lemma2}
\begin{proof} We can repeat the lines of the proof of lemma \ref{lemma:linindep_functions} with 
the only change at each selection step to choose $x_i \in U_i$. 

Since $f_1(x)$ can not be identically zero in $U_1$, there is $x_1\in U_1$ such that $f_1(x_1) \ne 0$. 
Assuming that there are $x_i\in U_i$,i=1,..,k-1, such that $det(A_1 := \{f_i(x_j)\}_{i=1,j=1}^{k-1,k-1}) \ne 0$, 
we can choose $x_k\in U_k$  such that for $x=x_k$ 
$$
det \left( \begin{array}{cccc}
f_1(x_1) & f_2(x_1) & ... & f_{k}(x_1) \\
... & ... & ... & ... \\
f_1(x_{k-1}) & f_2(x_{k-1}) & ... & f_{k}(x_{k-1}) \\
f_1(x) & f_2(x) & ... & f_{k}(x) \\
\end{array} \right) = \sum_{i=1}^{k-1} A_i f_i(x) \ne 0.
$$
Otherwise $\sum_{i=1}^{k} A_i f_i(x) = 0$, $\forall x\in U_k$ and $f_i$ would not be locally linear independent.
\qquad\end{proof} 
\newtheorem{rankmetric_ex2}[rankmetric_ex1]{Example}
\begin{rankmetric_ex2}
For any $k$ and distinct indices $s_1,...,s_k$, 
power functions $\{x^{s_l}\}_{l=1}^k$ are locally linear independent 
in any open $U\subset \mathbb{R}^n$.
Indeed, it is sufficient to show the claim for $n=1$ and $U=(a-\delta,a+\delta)$, $\delta>0$. 
Assume that $s_1<s_2<...<s_k$.
Let $\sum_{l=1}^k \alpha_l (x+a)^{s_l} = 0$, $\forall x\in(-\delta,\delta)$. 
Then the $(s_1)^{th}$, $(s_2)^{th}$,...,$(s_k)^{th}$ derivatives of the left-hand-side sum at $x=0$ have to be zeroes. 
Thus for any $m=1,...,k$, $$\sum_{l=m}^k \alpha_l \frac{s_l!}{(s_l-s_m)!}a^{s_l-s_m} = 0.$$
Which leads to the only possible choice $\alpha_l = 0$, for all $l$.
\end{rankmetric_ex2}
In fact, we have a stronger result.
\newtheorem{rankmetric_lemma3}[rankmetric_lemma1]{Lemma}
\begin{rankmetric_lemma3}
The power functions $\{x^s\}_{[s]\ge 0}^{\infty}$ are locally linear independent in $\mathbb{R}^n$.
\label{lemma:locally_linindep_powers}
\end{rankmetric_lemma3}
\begin{proof}
Let $U\subset \mathbb{R}^n$ be an open set.
Consider the Hilbert space $L_2(U)$ of square-integrable functions in $U$. 
We have that 
$$
L_2(U) \subset \overline{span\{x^s|_U,\textrm{ }[s] \ge 0\}},
$$
the closure of the span of the power functions (see \cite{riesz-nagy}, sec. 46). 
By the Gram-Schmidt method we can generate an 
orthonormal basis $\{\phi_s(x)\}_{[s]\ge 0}$ of $L_2(U)$ from $\{x^s|_U\}_{[s]\ge 0}$. 
Since $\{\phi_s(x)\}_{[s]\ge 0}$ can not be linear dependent, so can not be  $\{x^s|_U\}_{[s]\ge 0}$. 
\qquad\end{proof} 
%
%

For subsets $U_l$,l=1,...,m, of $\mathbb{R}^n$ with $\bigotimes_{l=1}^{m} U_l$ we denote 
the product $U_1 \times ... \times U_m$. 
\newtheorem{rankmetric_def3}[rankmetric_def1]{Definition}
\begin{rankmetric_def3}
Function $\psi:U\times V \to \mathbb{R}$ is said to have full rank almost everywhere ($a.e.$) in $U\times V$ if 
for any number $k\in\mathbb{N}$ and open sets $U_i\subset U$ and $V_j\subset V$, $i,j=1,...,k$ we have that 
$$
rank( \{\psi(x_i, y_j)\}_{i=1,j=1}^{k,k} ) \le k, 
$$
$\forall x_i\in U_i$ and $\forall y_j\in V_j$, 
with an equality for at least one choice of $x_i$'s and $y_j$'s.
\label{def:fullrank_function}
\end{rankmetric_def3}
\newtheorem{rankmetric_th2}[rankmetric_th1]{Proposition}
\begin{rankmetric_th2}
Let $\psi\in C^{\infty}(U\times V)$ be an analytic function in $V$ about a point $p\in V$ and 
for any $k\in\mathbb{N}$, there is a set of $k$ Taylor functions $c_s(x)$, which are locally linear independent in $U$.  
Then $\psi$ has full rank $a.e.$ in $U\times V$.
\label{theorem:rank_infinite}
\end{rankmetric_th2}
\begin{proof} Suppose that there exist $k$ and open sets $U_i\subset U$ and $V_j\subset V$, $i,j=1,...,k$, 
such that for all $x_i\in U_i$ and $y_j\in V_j$, $rank(\{\psi(x_i,y_j)\}_{i,j=1}^{k,k}) < k$.

If we fix $x_i\in U_i$, $i=1,...,k$, then by lemma \ref{lemma:locally_linindep_functions}, 
$\{\psi(x_i, y)\}_{i=1}^{k}$ are not locally linear independent functions of y in $V$. 
Consequently there exists an open set $W(x)\subset V$, $x=(x_1,...,x_{k})$, where $\{\psi(x_i, y)|_{W(x)}\}_{i=1}^{k}$ are 
linearly dependent. 
Thus there are $\alpha_i(x)$, not all zero, 
such that 
$$
\sum_{i=1}^{k} \alpha_i(x) \psi(x_i, y) = 0, \forall y\in W(x).
$$
If we expand $\psi(x_i,y)$'s into their Taylor series (assuming $p$ is the origin), we obtain 
$$
\sum_{s} \sum_{i=1}^{k} \alpha_i(x) c_s(x_i) y^s = 0, \forall y\in W(x).
$$
By local linear independence of the power functions $y^s$ it follows that 
$$
\sum_{i=1}^{k} \alpha_i(x) c_s(x_i) = 0, \forall s.
$$
Since the last condition is true for any $x_i\in U_i$, i=1,...,k, 
no $k$ functions $c_s$ would be locally linear indepent in $U$, which contradicts 
the theorem assumption. 
Therefore, $\psi$ has full rank in any open subset of $\bigotimes_{l=1}^{k} U \bigotimes_{l=1}^{k} V$. 
\qquad\end{proof} 

The following fact will be useful.
\newtheorem{rankmetric_lemma4}[rankmetric_lemma1]{Lemma}
\begin{rankmetric_lemma4}
Any collection $\mathcal{O}$ of (nonempty) open and disjoin subsets of $\mathbb{R}^n$ is countable.
\label{lemma:union_open_sets}
\end{rankmetric_lemma4}
{\it Proof.}
Let $O\in\mathcal{O}$. Then there is at least one vector $q\in O$ with all rational coordinates. 
By applying the axiom of choice we define a map 
$
O \mapsto (k_1,l_1,k_2,l_2,...,k_n,l_n),
$
where $k_i$, $l_i$, $i=1,...,n$, are whole numbers such that 
$q=(2^{k_1}l_1,...,2^{k_n}l_n) \in O$. 
$\mathcal{O}$ is countable since this map is injective one 
from $\mathcal{O}$ to $\mathbb{Q}^{2n}$.
\qquad$\Box$

Let $\nu_n$ be the Lebesque measure in $\mathbb{R}^n$. 
The next measure-theoretic result has important consequences.
\newtheorem{rankmetric_th13}[rankmetric_th1]{Proposition}
\begin{rankmetric_th13}
If $f\in C^{\infty}(U)$ for an open subset $U$ of $\mathbb{R}^n$, 
such that the set $V_0 := \{x\in U| f(x) = 0\}$ has no interior points, then $\nu_n(V_0) = 0$.
\label{theorema:f_nullset_zeros}
\end{rankmetric_th13}
{\it Proof.}
Define $V_1 := \{x\in U| \frac{\partial f}{\partial x}|_{x} = 0\}$.
When $n>1$, $\frac{\partial f}{\partial x}$ is a gradient vector and $\frac{\partial f}{\partial x}|_{x} = 0$ 
means that all partial derivatives at point $x$ vanish. 
Both $V_0$ and $V_1$ are closed in $U$. 
By the implicit function theorem, for any point $x\in V_0\backslash V_1$ there is 
an open ball $B_x=B_x(r_x)$ of radius $r_x>0$ and 
a hypersurface $S_x$ of dimension $n-1$, the graph of a function of $n-1$ variables, 
such that $x\in S_x = B_x \cap V_0\backslash V_1$, 
Therefore there exists a collection $\{S_{\alpha}\}_{\alpha}$ that partitions $V_0\backslash V_1$, 
$V_0\backslash V_1 = \cup_{\alpha} S_{\alpha}$, and such that 
$S_{\alpha}$ are disjoin and with measure zero, $\nu_n(S_{\alpha}) = 0$. 
Each $S_{\alpha}$ is obtained by merging countable number of hypersurfaces $S_x$.
Indeed, let for every $x$, $U_x$ be the support of the function which graph is $S_x$. 
We may assume that every $U_x$ is a subset of one of the coordinate hyperplane $x^i = 0$, $i=1,...,n$. 
Those of them that are subsets of $x^i = 0$ form a countable disjoin collection of open subsets, 
since any two $U_x$ on $x^i = 0$ that intersect with each other 
can be seemlessly merged into the support of one function. 

Observe that if $x\in S_{\alpha}$ and $y\in S_{\beta}$ for $\alpha\ne\beta$ then 
$B_x(r_x/2)\cap B_y(r_y/2) = \emptyset$, 
for otherwise either $x\in B_y(r_y)$ or $y\in B_x(r_x)$, a contradiction. 
Therefore for the open sets $U_{\alpha} := \cup_{x\in S_{\alpha}} B_x(r_x/2)$, 
we have $U_{\alpha} \cap V_0\backslash V_1 = S_{\alpha}$ and $U_{\alpha}$ are disjoin. 
By the virtue of Lemma \ref{lemma:union_open_sets} the collection of open sets 
$\{U_{\alpha}\}_{\alpha}$ is countable and so is $\{S_{\alpha}\}_{\alpha}$. 
Therefore $\nu_n(V_0\backslash V_1) = 0$ and $\nu_n(V_0) = \nu_n(V_1\cap V_0)$.

We may repeat the same analysis for higher derivatives of order $l > 1$ 
and apply induction on $l$. Define 
$$
V_l := \{x\in U | \textrm{For all }s=(s_1,...,s_n),\textrm{ s.t. }[s] = l:\textrm{ } \frac{\partial^s f}{\partial x^s}|_{x} = 0\}.
$$
For example, let $l=2$ and let 
$$
V_1^i := \{x\in U| \frac{\partial f}{\partial x^i}|_{x} = 0\}, ~  
V_2^i := \{x\in U| \frac{\partial^2 f}{\partial x^i\partial x^j}|_{x} = 0,\textrm{ for all }j=1,...,n\}.
$$
The above argument for showing $\nu_n(V_0) = \nu_n(V_1\cap V_0)$ is applicable 
to all pairs $(V_1^i, V_2^i)$ and consequently $\nu_n(V_1^i) = \nu_n(V_2^i\cap V_1^i)$, $i=1,...,n$. 
Next we use the fact that for any four measurable sets satisfying $B\subset A$, $\nu_n(B)=\nu_n(A)$, 
$D\subset C$ and $\nu_n(D)=\nu_n(C)$ we have $\nu_n(B\cap D)=\nu_n(A\cap C)$, 
to conclude that 
$\nu_n(\cap_{i=1}^n V_1^i) = \nu_n(\cap_{i=1}^n V_1^i \cap_{i=1}^n V_2^i)$ or 
equivalently $\nu_n(V_1) = \nu_n(V_1 \cap V_2)$.

By induction, $\nu_n(V_{l-1}\backslash V_l) = 0$ for any $l>0$. 
Therefore given that $\nu_n(V_0) = \nu_n(\cap_{m=0}^{l-1} V_m)$, 
from 
$$
\nu_n(\cap_{m=0}^{l-1} V_m) = \nu_n((\cap_{m=0}^{l-1} V_m) \backslash V_l) + \nu_n(\cap_{m=0}^{l} V_m) = \nu_n(\cap_{m=0}^{l} V_m) 
$$ 
we obtain $\nu_n(V_0) = \nu_n(\cap_{m=0}^{l} V_m)$.
Finally, we realize that the set $\cap_{s=1}^{\infty} V_s$ should be empty, 
for otherwise, with all vanishing derivatives at one point, the function $f$ would vanish in an open subset of $U$.
Thus we conclude that $\nu_n(V_0) = \nu_n(\cap_{l=1}^{\infty} V_l) = 0$.
\qquad$\Box$

By applying the above proposition for the function $det(\{\psi(x_i,y_j)\}_{i,j=1}^{k,k})$ we immediately obtain the following result.
\newtheorem{rankmetric_cor1}{Corollary}
\begin{rankmetric_cor1}
Under the conditions of Theorem \ref{theorem:rank_infinite}, for any $k\in\mathbb{N}$ 
$$
\nu_{2kn}(D_k := \{(x_1,...x_k,y_1,...,y_k)\in \bigotimes_{l=1}^{k} U \bigotimes_{l=1}^{k} V: rank(\{\psi(x_i,y_j)\}_{i,j=1}^{k,k}) < k\}) = 0, 
$$
which explains "$a.e.$" term in the full rank definition.
\label{corollary:rank_infinite}
\end{rankmetric_cor1}

Next we show that Theorem \ref{theorem:rank_infinite} can be applied to a large class 
of symmetric analytic functions. 
With $x.y = \sum_{l=1}^n x^l y^l$ we will denote the dot product in $\mathbb{R}^n$.
\newtheorem{rankmetric_cor2}[rankmetric_cor1]{Corollary}
\begin{rankmetric_cor2}
Let $\psi(x,y) = h(x.y)$, $(x,y)\in U\times V$ 
and for $W=\{x.y|x\in U, y\in V\}\subset \mathbb{R}$, $0\in W$.
Let also $h$ be analytic function in $W$ 
about the origin $0$ and such that $h^{(s)}(0)\ne 0$ for infinitely many $s$.
Then $\psi$ has full rank $a.e.$ in $U \times V$. 
\label{corollary:symm_rank_infinite}
\end{rankmetric_cor2}
\begin{proof}
Observe that 
$$
c_s(x) = \frac{\partial^s}{\partial y^s} \psi(x,y)|_{y=0} = x^s h^{(s)}(0).
$$
Since power functions are locally linear independent, for any $k$ there exist 
$k$ functions $c_s$ that are locally linear independent in $U$. 
Thus the conditions of Theorem \ref{theorem:rank_infinite} are fulfilled.
\qquad\end{proof} 

\newtheorem{rankmetric_ex3}[rankmetric_ex1]{Example}
\begin{rankmetric_ex3}
Consider the distance on the n-sphere $\mathbb{S}^n$ given by 
$$
\psi(x,y) = \cos^{-1}(x.y), \\ x,y\in \mathbb{S}^n.
$$
Conditions of corollary \ref{corollary:symm_rank_infinite} are met for 
the function $h(z) = \cos^{-1}(z)$, $z\in[-1,1]$. 
Indeed $h$ is analytic about $0$ with infinitely many Taylor coefficients non-zero
$$
\cos^{-1}(z) = \frac{\pi}{2} - \sum_{k=0}^{\infty} \frac{(2k)!}{2^{2k}(k!)^2} \frac{z^{2k+1}}{2k+1}, \textrm{ } |z|<1.
$$
\label{ex:full_rank_sphere}
\end{rankmetric_ex3}

Example \ref{ex:full_rank_sphere} shows that the distance on the n-sphere is 
of full rank almost everywhere, a fact of importance to some inverse problems on the sphere.

\newtheorem{rankmetric_remark2}[rankmetric_remark1]{Remark}
\begin{rankmetric_remark2}
The full rank property of $\cos^{-1}(x.y), \textrm{ } x,y\in \mathbb{S}^n$ follows from another standard result in 
linear operator theory. In the Hilbert space $L_2(\mathbb{S}^n)$ of square integrable functions on $\mathbb{S}^n$, 
the operator with kernel $\cos^{-1}(x.y)$ is symmetric and thus we have the representation $\cos^{-1}(x.y) = \sum_{k\ge 1}\lambda_k\phi_k(x)\phi_k(y)$, 
where $\lambda_k$ and $\phi_k$ are the eigenvalues and eigenvectors of this operator (see Theorems of Hilbert and Schmidt, \cite{riesz-nagy}, sec. 97). 
Since $\phi_k$'s form orthonormal system in $L_2(\mathbb{S}^n)$, they are necessarily linear independent. 
In fact, one can expect that they are locally linear independent, which will give us the full rank result. 
Our approach, however, is stronger in a sense. 
For example, according to corollary \ref{corollary:symm_rank_infinite}, $cos(x.y), \textrm{ } x,y\in \mathbb{R}^n$, has full rank $a.e.$, but 
one can not use Hilbert space argument to show it, since $cos(x.y)$ is not square integrable in $\mathbb{R}^n\times\mathbb{R}^n$.
\end{rankmetric_remark2}

\section{On solvability of some systems of matrix equations}

Let $U$ and $V$ be open subsets of $\mathbb{R}^n$ and 
$\eta:U\times U \to \mathbb{R}^n$ be a continuous vector valued function. 
Let also $\mathcal{P}=\{p_i\in U\}_{i=1}^k$ be a set of points in $U$.
For any discrete function $f$ on $\mathcal{P}$, $f=(f_i=f(x_i))_{i=1}^k \in\mathbb{R}^k$, define $\Sigma[f]=\{\Sigma[f]_j\}_{j=1}^k$, 
$$
\Sigma[f]_j := \sum_{i=1}^k f_i \eta(x_j,x_i)\eta(x_j,x_i)^T,\textrm{ }j=1,...,k.
$$ 
We have $\Sigma[f]\in [Sym_n]^k$, where $Sym_n$ is the set of symmetric $n\times n$ matrices.
Note that $K(x,y) := \eta(x,y)\eta(x,y)^T$ is always symmetric and positive semi-definite and 
so $K$ is Mercer kernel (see \cite{cucker-smale-regul}). 

If $C = \Sigma[f^0]$ for $f^0 \in\mathbb{R}^n$, then $f^0$ will be a solution  of the system 
\begin{equation}\label{eq:rankmetric_f0_equation}
\sum_{i=1}^k f_i Y_{ji} = C_j,\\ j=1...k,
\end{equation}
where $Y_{ji} = \eta(x_j,x_i)\eta(x_j,x_i)^T \in Sym_n$.
We say that $f^0$ can be recovered if system (\ref{eq:rankmetric_f0_equation}) 
has a unique solution. 
System (\ref{eq:rankmetric_f0_equation}) is in fact a discretization of the integral tensor equation (\ref{eq:rankmetric_inverse_tensor_field}).

We re-arrange the components of each of the $k^2$ $n\times n$-matrices $Y_{ji} = \{Y_{ji}^{lm}\}_{l=1,m=1}^{n,n}$, i,j=1,...,k , 
into new $k\times k$-matrices $B_{lm}=\{Y_{ji}^{lm}\}_{i=1,j=1}^{k,k} \in \mathbb{R}^{k\times k}$, l,m,=1,...,n.
Then we order all $n^2$ matrices $B_{lm}$ vertically to obtain a $n^2k\times k$-matrix
$$
\mathcal{Y} = [B_{11}^T ... B_{1n}^T .... B_{n1}^T...B_{nn}^T]^T \in \mathbb{R}^{n^2k \times k}.
$$
Similarly, matrices $C_j$ can be unfolded in a large $n^2k$-vector
$$
\mathcal{C} = [C_1^{11} ... C_k^{11} C_1^{12} ... C_k^{12}... ... C_1^{nn} ... C_k^{nn}]^T, 
$$
and equation (\ref{eq:rankmetric_f0_equation}) can be written as 
$$
\mathcal{Y}f = \mathcal{C}, \\ f\in\mathbb{R}^n.
$$
Occasionally we will write $\mathcal{Y}(\mathcal{P})$ to emphasize that $\mathcal{Y}$ 
is contingent on the choice of $\mathcal{P}$.

In fact, when $C = \Sigma[f^0]$, $f^0$ can be recovered if and only if 
\begin{equation}\label{eq:f0_equation2}
\sum_{i=1}^k f_i tr(Y_{ji}) = tr(C_j),\textrm{ } j=1...k, 
\end{equation}
has a unique solution ($f^0$). Define 
$$
\psi(x,y) := tr(\eta(x,y)\eta(x,y)^T) = ||\eta(x,y)||^2
$$ and let 
$\Psi = \{\psi(x_i,y_j)\}_{i,j=1}^{k,k}$ and $c=(tr(C_1),...,tr(C_k))^T$.
Then (\ref{eq:f0_equation2}) can be written as $\Psi f = c$.
The rank of function $\psi$, which we studied in the previous section, determines the solvability of (\ref{eq:f0_equation2}).

Since $C_j=\Sigma[f^0]_j \in span(Y_{j1},...,Y_{jk})$, we have that 
$$
rank(\mathcal{Y}|\mathcal{C}) = rank(\mathcal{Y}), 
$$
where $\mathcal{Y}|\mathcal{C}$ is matrix 
$\mathcal{Y}$ with vector $\mathcal{C}$ attached as a last column.
The system of linear equations (\ref{eq:rankmetric_f0_equation}) has 
a unique solution if and only if $\mathcal{Y}$ is of full rank, $rank(\mathcal{Y}) = k$.

Therefore it is of importance to find under what circumstances $\mathcal{Y}$ has full rank.
If the points of $\mathcal{P}$ are chosen by a continuous distribution of M, 
what will be the probability for $\mathcal{Y}$ to be of smaller rank?

The properties of $\mathcal{Y}$ are contingent on the choice of the vector function $\eta$.
It is natural first to consider the simplest choice $\eta(x,y) = x-y$, corresponding to the Euclidean distance vector function.
Observe that in this case $\psi(x,y) = tr((y-x)(y-x)^T) = ||y-x||^2$ and 
as we already have found in the previous section (Example \ref{ex:rank_Rn}), 
$rank(\psi) = n+2$. Thus 
$$
rank(\Psi) \le n+2.
$$
In the light of this fact, we may expect that the rank of $\mathcal{Y}$ is also bounded for this choice of $\eta$. 
\newtheorem{rankmetric_th3}[rankmetric_th1]{Proposition}
\begin{rankmetric_th3}
For the map $\eta(x,y) = x-y$ in $\mathbb{R}^n$, 
$$
rank(\mathcal{Y}(\mathcal{P})) = \min\{k,\frac{(n+1)(n+2)}{2}\},
$$
and therefore, (\ref{eq:rankmetric_f0_equation}) has not a unique solution for $k > \frac{(n+1)(n+2)}{2}$.
\label{theorem:rank_Rn}
\end{rankmetric_th3}
\begin{proof}
With respect to a fixed coordinate system in $\mathbb{R}^n$, 
we identify the points $p_i\in\mathcal{P}$ with vectors $v_i = (v_i^1,...,v_i^n)\in\mathbb{R}^n$.
Then the matrices $Y_{ji}$ have the following structure 
\begin{equation}\label{def:matrix_Yij}
Y_{ji} = \left( \begin{array}{cccc}
(v_i^1-v_j^1)^2 & (v_i^1-v_j^1)(v_i^2-v_j^2) & ... & (v_i^1-v_j^1)(v_i^n-v_j^n) \\
(v_i^2-v_j^2)(v_i^1-v_j^1) & (v_i^1-v_j^1)^2 & ... & (v_i^2-v_j^2)(v_i^n-v_j^n) \\
... & ... & ... & ... \\
(v_i^n-v_j^n)(v_i^1-v_j^1) & (v_i^n-v_j^n)(v_i^2-v_j^2) & ... & (v_i^n-v_j^n)^2
\end{array} \right)
\end{equation}

First we look at the case $n=1$, because is much simpler and gives us intuition needed for the general case. 
We will show that for any eight real numbers $a,b,c,d,a_1,b_1,c_1,d_1$, 
the matrix 
\begin{equation}\label{def:matrix_abcd}
\mathbf{X} = \left( \begin{array}{cccc}
(a-a_1)^2 & (b-a_1)^2 & (c-a_1)^2 & (d-a_1)^2 \\
(a-b_1)^2 & (b-b_1)^2 & (c-b_1)^2 & (d-b_1)^2 \\
(a-c_1)^2 & (b-c_1)^2 & (c-c_1)^2 & (d-c_1)^2 \\
(a-d_1)^2 & (b-d_1)^2 & (c-d_1)^2 & (d-d_1)^2
\end{array} \right)
\end{equation}
is singular with $rank(\mathbf{X}) \le 3$. 
Since any $4\times 4$ sub-matrix of $\mathcal{Y}$ has the structure of $\mathbf{X}$, 
it will follow that $rank(\mathcal{Y}) \le 3$.

For any $k\times k$ matrix $Z$ let $l^j(Z)$ denote the j-th column vector of $Z$ 
and $\mathcal{L}(Z)=span\{l^1(Z),...,l^k(Z)\}$, the linear space defined by the column vectors of $Z$.

Too see why $rank(\mathbf{X}) \le 3$, 
we define $\underline{1}=(1,1,1,1)^T$, $\underline{x}=(a_1,b_1,c_1,d_1)^T$,
$\underline{x}^2=(a_1^2,b_1^2,c_1^2,d_1^2)^T$, 
and consider following expression for the columns of $X$ 
$$
l^1(X) = a^2\underline{1} - 2a\underline{x} + \underline{x}^2, \textrm{ }
l^2(X) = b^2\underline{1} - 2b\underline{x} + \underline{x}^2 
$$
$$
l^3(X) = c^2\underline{1} - 2d\underline{x} + \underline{x}^2, \textrm{ } 
l^4(X) = d^2\underline{1} - 2c\underline{x} + \underline{x}^2 .
$$
Then we observe that 
$
\mathcal{L}(X) \subset \mathcal{X} := span\{\underline{1},\underline{x},\underline{x}^2\}, 
$
and therefore $rank(X) \le dim(\mathcal{X}) \le 3$.

Let check now the general case. 
The $n^2k\times k$ matrix $\mathcal{Y}$ looks like this 
$$
\mathcal{Y} = \left( \begin{array}{cccc}
(v_1^1-v_1^1)(v_1^1-v_1^1) & (v_1^1-v_2^1)(v_1^1-v_2^1) & ... & (v_1^1-v_k^1)(v_1^1-v_k^1) \\
(v_2^1-v_1^1)(v_2^1-v_1^1) & (v_2^1-v_2^1)(v_2^1-v_2^1) & ... & (v_2^1-v_k^1)(v_2^1-v_k^1) \\
... & ... & ... & ... \\
(v_k^1-v_1^1)(v_k^1-v_1^1) & (v_k^1-v_2^1)(v_k^1-v_2^1) & ... & (v_k^1-v_k^1)(v_k^1-v_k^1) \\

(v_1^1-v_1^1)(v_1^2-v_1^2) & (v_1^1-v_2^1)(v_1^2-v_2^2) & ... & (v_1^1-v_k^1)(v_1^2-v_k^2) \\
(v_2^1-v_1^1)(v_2^2-v_1^2) & (v_2^1-v_2^1)(v_2^2-v_2^2) & ... & (v_2^1-v_k^1)(v_2^2-v_k^2) \\
... & ... & ... & ... \\
(v_k^1-v_1^1)(v_k^2-v_1^2) & (v_k^1-v_2^1)(v_k^2-v_2^2) & ... & (v_k^1-v_k^1)(v_k^2-v_k^2) \\ 
... & ... & ... & ... \\
(v_1^1-v_1^1)(v_1^n-v_1^n) & (v_1^1-v_2^1)(v_1^n-v_2^n) & ... & (v_1^1-v_k^1)(v_1^n-v_k^n) \\
(v_2^1-v_1^1)(v_2^n-v_1^n) & (v_2^1-v_2^1)(v_2^n-v_2^n) & ... & (v_2^1-v_k^1)(v_2^n-v_k^n) \\
... & ... & ... & ... \\
(v_k^1-v_1^1)(v_k^n-v_1^n) & (v_k^1-v_2^1)(v_k^n-v_2^n) & ... & (v_k^1-v_k^1)(v_k^n-v_k^n) \\ 
... & ... & ... & ... \\

... & ... & ... & ... \\
(v_1^n-v_1^n)(v_1^1-v_1^1) & (v_1^n-v_2^n)(v_1^1-v_2^1) & ... & (v_1^n-v_k^n)(v_1^1-v_k^1) \\
(v_2^n-v_1^n)(v_2^1-v_1^1) & (v_2^n-v_2^n)(v_2^1-v_2^1) & ... & (v_2^n-v_k^n)(v_2^1-v_k^1) \\
... & ... & ... & ... \\
(v_k^n-v_1^n)(v_k^1-v_1^1) & (v_k^n-v_2^n)(v_k^1-v_2^1) & ... & (v_k^n-v_k^n)(v_k^1-v_k^1) \\ 
... & ... & ... & ... \\
(v_1^n-v_1^n)(v_1^n-v_1^n) & (v_1^n-v_2^n)(v_1^n-v_2^n) & ... & (v_1^n-v_k^n)(v_1^n-v_k^n) \\
(v_2^n-v_1^n)(v_2^n-v_1^n) & (v_2^n-v_2^n)(v_2^n-v_2^n) & ... & (v_2^n-v_k^n)(v_2^n-v_k^n) \\
... & ... & ... & ... \\
(v_k^n-v_1^n)(v_k^n-v_1^n) & (v_k^n-v_2^n)(v_k^n-v_2^n) & ... & (v_k^n-v_k^n)(v_k^n-v_k^n) \\ 

\end{array} \right)
$$

Consider the $s$-th column of $\mathcal{Y}$. Its elements are formed by multiplication of 
two differences, which we will express formally as differences of sub-columns $(1)-(2)$ and $(3)-(4)$.
Next we show the four sub-columns of $s$-th column as rows 
$$
(1)^T: \underbrace{\underbrace{v_1^1,...v_k^1}_{k}...\underbrace{v_1^1,...v_k^1}_{k}}_{n} 
\underbrace{\underbrace{v_1^2,...v_k^2}_{k}...\underbrace{v_1^2,...v_k^2}_{k}}_{n} ...
\underbrace{\underbrace{v_1^n,...v_k^n}_{k}...\underbrace{v_1^n,...v_k^n}_{k}}_{n}
$$
$$
(2)^T: \underbrace{\underbrace{v_s^1,...v_s^1}_{k}...\underbrace{v_s^1,...v_s^1}_{k}}_{n} 
\underbrace{\underbrace{v_s^2,...v_s^2}_{k}...\underbrace{v_s^2,...v_s^2}_{k}}_{n} ...
\underbrace{\underbrace{v_s^n,...v_s^n}_{k}...\underbrace{v_s^n,...v_s^n}_{k}}_{n}
$$
$$
(3)^T: \underbrace{v_1^1,...v_k^1}_{k}...\underbrace{v_1^n,...v_k^n}_{k}
\underbrace{v_1^1,...v_k^1}_{k}...\underbrace{v_1^n,...v_k^n}_{k} ...
\underbrace{v_1^1,...v_k^1}_{k}...\underbrace{v_1^n,...v_k^n}_{k}
$$
$$
(4)^T: \underbrace{v_s^1,...v_s^1}_{k}...\underbrace{v_s^n,...v_s^n}_{k}
\underbrace{v_s^1,...v_s^1}_{k}...\underbrace{v_s^n,...v_s^n}_{k} ...
\underbrace{v_s^1,...v_s^1}_{k}...\underbrace{v_s^n,...v_s^n}_{k}
$$
Formally, $s$-th column equals $((1)-(2))((3)-(4))$ if the operations are taken component-wise.

For a vector $x$, with $x\hookrightarrow_m$ we denote the vector obtained from $x$ by 
shifting m positions in the right.
For example, $(0,1,0,0)\hookrightarrow_2 = (0,0,0,1)$.
Let 
$$
\underline{x}_1 = (\underbrace{v_1^1,...v_k^1,v_1^2,...v_k^2,...,v_1^n,...v_k^n}_{nk},\underbrace{0,...,0}_{(n-1)nk})^T \in \mathbb{R}^{n^2k},
$$
$\underline{x}_j= \underline{x}_1 \hookrightarrow_{(j-1)nk}$, for j=1,...n,
$$
\underline{z}_1 = (v_1^1,...v_k^1\underbrace{0,...,0}_{(n-1)k}v_1^2,...v_k^2\underbrace{0,...,0}_{(n-1)k},...,v_1^n,...v_k^n\underbrace{0,...,0}_{(n-1)k})^T \in \mathbb{R}^{n^2k},
$$
$\underline{z}_j = \underline{z}_1 \hookrightarrow_{(j-1)k}$, for j=1,...n,
$$
\underline{1}_{1} = (\underbrace{1,...,1}_{k},\underbrace{0,...,0}_{(n-1)k},\underbrace{0,...,0}_{(n-1)nk})^T \in \mathbb{R}^{n^2k},
$$
$\underline{1}_{j} = \underline{1}_{1} \hookrightarrow_{(j-1)}$, for j=1,...n, and 
$$
\underline{1}_{ij} = 
(\underline{1}_1 \hookrightarrow_{(i-1)nk}) + (\underline{1}_1 \hookrightarrow_{(i-1)nk+(j-1)k}),
\textrm{ for }j=1,...,n,\\ i<j.
$$
Let $l^0\in\mathbb{R}^{n^2k}$ be the vector
$$
(v_1^1v_1^1,v_2^1v_2^1,...,v_k^1v_k^1,v_1^1v_1^2,v_2^1v_2^2,...,v_k^1v_k^2,...,v_1^1v_1^n,v_2^1v_2^n,...,v_k^1v_k^n,
$$
$$
v_1^2v_1^1,v_2^2v_2^1,...,v_k^2v_k^1,v_1^2v_1^2,v_2^2v_2^2,...,v_k^2v_k^2,...,v_1^2v_1^n,v_2^2v_2^n,...,v_k^2v_k^n, ... 
$$
$$
v_1^nv_1^1,v_2^nv_2^1,...,v_k^nv_k^1,v_1^nv_1^2,v_2^nv_2^2,...,v_k^nv_k^2,...,v_1^nv_1^n,v_2^nv_2^n,...,v_k^nv_k^n)^T.
$$
Then we can express $s$-th column of $\mathcal{Y}$ as
$$
l^s(\mathcal{Y}) = 
l^0
- \sum_{j=1}^n v_s^j (\underline{x}_j + \underline{z}_j)
+ \sum_{j=1}^n  (v_s^j)^2 \underline{1}_j + 2\sum_{i<j}v_s^iv_s^j \underline{1}_{ij}.
$$
Note that $\underline{1}_j = \underline{1}_{1j}$ and therefore 
$$
\mathcal{L}(\mathcal{Y}) \subset span\{l^0, \underline{1}_{ij},(\underline{x}_j+\underline{z}_j)\}_{j=1,i<j}^n.
$$
Finally we obtain 
$$
dim(\mathcal{L}(\mathcal{Y})) \le 1 + \frac{n(n+1)}{2} + n = \frac{(n+1)(n+2)}{2}.
$$
\qquad\end{proof} 

Another matrix of potential interest is obtained by arranging 
matrices $Y_{ji}$,i,j=1,...,k in one $kn\times kn$ matrix 
\begin{equation}\label{def:matrix_Z}
\mathcal{Z} = \left( \begin{array}{cccc}
Y_{11} & Y_{21} & ... & Y_{k1} \\
Y_{12} & Y_{22} & ... & Y_{k2} \\
... & ... & ... & ... \\
Y_{1k} & Y_{21} & ... & Y_{kk} \\
\end{array} \right)
\end{equation}
What is the rank of $\mathcal{Z}$ for the choice $\eta(x,y) = y-x$? 
\newtheorem{rankmetric_th4}[rankmetric_th1]{Proposition}
\begin{rankmetric_th4}
For the map $\eta(x,y) = x-y$ in $\mathbb{R}^n$, 
$$
rank(\mathcal{Z}) = \min\{k,n(n+2)\}.
$$
\label{hypothesis:rank_Z_Rn}
\end{rankmetric_th4}
\begin{proof}
We will find a representation of $(s-1)n+m$-th column of $\mathcal{Z}$ for s=1,...,k and m=1,...,n,
as a linear combination of $n(n+2)$ basis column vectors and that will proof the claim.

The elements of $(s-1)n+m$-th column of $\mathcal{Y}$ are formed by multiplication of 
two differences, which we express formally as differences of sub-columns $(1)-(2)$ and $(3)-(4)$.
We show these four sub-columns as rows 
$$
(1)^T: \underbrace{\underbrace{v_s^1,...v_s^n}_{n}\underbrace{v_s^1,...v_s^n}_{n}...\underbrace{v_s^1,...v_s^n}_{n}}_{k} 
,\textrm{   }
(2)^T: \underbrace{v_1^1,...v_1^n}_{n}\underbrace{v_2^1,...v_2^n}_{n}...\underbrace{v_k^1,...v_k^n}_{n} 
$$
$$
(3)^T: \underbrace{\underbrace{v_s^l,...v_s^l}_{n}\underbrace{v_s^l,...v_s^l}_{n}...\underbrace{v_s^l,...v_s^l}_{n}}_{k} 
,\textrm{   }
(4)^T: \underbrace{v_1^l,...v_1^l}_{n}\underbrace{v_2^l,...v_2^l}_{n}...\underbrace{v_k^l,...v_k^l}_{n} 
$$
Then, $(s-1)n+m$-th column equals $((1)-(2))((3)-(4))$ if the operations are taken component-wise.
Let
$$
\underline{x}_1^m = (\underbrace{v_1^m,0...0}_{n},\underbrace{v_2^m,0...0}_{n},...,\underbrace{v_k^m,0...0}_{n})^T \in \mathbb{R}^{nk},
$$
$\underline{x}_j^m = \underline{x}_1^m \hookrightarrow_{(j-1)}$, for j=1,...n,
$$
\underline{w} = (\underbrace{v_1^1,...,v_1^n}_{n}\underbrace{v_2^1,...,v_2^n}_{n},...,\underbrace{v_k^1,...,v_k^n}_{n} )^T \in \mathbb{R}^{nk},
$$
$$
\underline{1}_{1} = (\underbrace{1,0...0}_{n},\underbrace{1,0...0}_{n},...,\underbrace{1,0...0}_{n})^T \in \mathbb{R}^{nk},
$$
$\underline{1}_{j} = \underline{1}_{1} \hookrightarrow_{(j-1)}$, for j=1,...n,
and
$$
\underline{l}^m = (v_1^1v_1^m,...,v_1^n v_1^m,  v_2^1v_2^m,...,v_2^nv_2^m,...,v_k^1v_k^m,...,v_k^nv_k^m)^T \in\mathbb{R}^{nk}
$$
Now we can write $(s-1)n+m$-th column of $\mathcal{Z}$ as
$$
l^{(s-1)n+m}(\mathcal{Z}) = 
\underline{l}^m
+ \sum_{j=1}^n v_s^j \underline{x}_j^m
- v_s^m \underline{w} 
- v_s^m \sum_{j=1}^n v_s^j \underline{1}_{j}.
$$
Observe that $\underline{w} = \sum_{j=1}^n \underline{x}_j^j$.
Therefore 
$$
\mathcal{L}(\mathcal{Z}) \subset span\{\underline{l}^m, \underline{1}_{j}, \underline{x}_j^m\}_{j=1,m=1}^{n,n}, 
$$
and
$$
dim(\mathcal{L}(\mathcal{Z})) \le n + n + n^2 = n(n+2).
$$
\qquad\end{proof} 

In conclusion, the rank of matrices $\mathcal{Y}$ and $\mathcal{Z}$ depend on the rank 
of the bi-variate function $\psi$. If $\psi$ has a finite rank then the ranks of 
$\mathcal{Y}$ and $\mathcal{Z}$ will be bounded as is the case with the 
Euclidean distance function $\psi(x,y) = ||x-y||^2$.
As a consequence, a distribution function $f$ in $\mathbb{R}^n$ can not be recovered from its covariance field $\Sigma[f]$.

Recalling Corollary \ref{corollary:rank_infinite} we immediately obtain the following 
\newtheorem{rankmetric_cor3}[rankmetric_cor1]{Corollary}
\begin{rankmetric_cor3}
If $\psi(x,y)=tr(\eta(x,y)\eta(x,y)^T)$ has full rank $a.e.$ in $U$ 
and $\mathcal{P}$ is a finite sample of points drawn by a continuous distribution on $U$, then 
$\mathcal{Y(P)}$ has full rank with probability one.
\label{corollary:Y_rank_infinite}
\end{rankmetric_cor3}

According to Example \ref{ex:full_rank_sphere}, the square sphere distance, 
$\psi(x,y) = (\cos^{-1}(x.y))^2$, is of full rank and Corollary \ref{corollary:Y_rank_infinite} holds.
Thus, in contrast to the Euclidean space, 
a function $f$ on $\mathbb{S}^n$ can always be recovered from the respected field $\Sigma[f]$.

\section{Final discussion and simulation results}

Here we briefly discuss the problem of regularization of ill-posed integral equations 
of type (\ref{eq:rankmetric_inverse_trace_field}). We re-formulate it using the geodesic distance $d$ 
on a Riemannian manifold with volume form $V$
$$
\int_U d^2(p,q) f(q) dV(q) = g(p),\textrm{ }U\subset M.
$$
Discretizing the above equation on a set $\{p_i\}_i$ of $k$ points results to a system of linear equations $D\tilde f = \tilde g$, 
where $D=\{d^2(p_i,p_j)\}_{i,j=1}^{k,k}$. Although for a full-rank metric $d$ on M, 
$D$ will be almost always non-singular, its 
condition number, $cond(D)=||D||\textrm{ }||D^{-1}||$ , may increase dramatically as $k$ increases, 
making solutions unstable. 

Regularization of ill-posed linear problems is an area of active research (see \cite{ivanov}, \cite{tikhonov}, \cite{cucker-smale-regul} and \cite{cavalier-tsybakov}). 
Along with the classical Tikhonov regularization several new approaches are currently investigated. 
A good overview of the topic can be found in \cite{cavalier-inverse}.

Instead of trying to regularize an ill-conditioned kernel $\psi=d^2$, 
we replace it with a well-conditioned element from a family of similar to it kernels.
The procedure is applicable to all non-negative kernels. 

Consider the class of kernels $\{(d(p,q)-\alpha)^2\}_{\alpha>0}$. 
What would be a good element of this class that gives low (in average) conditional number for $D$? 
Assuming $U$ to be compact, we choose $\alpha$ that minimizes the expectation of 
$(d(X,Y)-\alpha)^2$ for $X,Y\sim Unif(U)$, and $X$ and $Y$ independent. 
It is a well known statistical fact that 
the optimal choice in this case is the expectation of $d(X,Y)$.

For example, on the unit sphere $\mathbb{S}^n$ if we take $U$ to be the whole manifold, 
then $E(d(X,Y)) = \pi/2$ and 
therefore, $\alpha = \pi/2$ is the optimal choice for it minimizes $E(d(X,Y)-\alpha)^2$.
Next we show some empirical evidence confirming the efficiency of this choice.

For a series of values of $k$ ranging from 50 to 900 we draw 20, $k$-samples on $\mathbb{S}^2$, 
calculate their conditional numbers, 
and then report the mean conditional number for these 20 samples. The next table shows how the mean conditional number varies 
as $k$ increases for the baseline choice $\alpha=0$ and the alternative $\alpha=\pi/2$.
\begin{center}
\begin{tabular}{|l|l|l|l|l|l|l|l|}
\hline
$\alpha$ & k=50 & k=100 & k=150 & k=200 & k=250 & k=300 &  \\ \hline
$0$     & 0.0000 & 0.0002 & 0.0226 & 0.3745 & 3.4871s & 1.1612 & $\times 10^{+15}$ \\ 
$\pi/2$ & 0.1291 & 0.0710 & 0.2865 & 4.1476 & 0.5755 & 0.1442 & $\times 10^{+4}$\\ \hline
        & k=400 & k=500 & k=600 & k=700 & k=800 & k=900 &  \\ \hline
$0$     & 0.0264 & 0.0694 & 0.0562 & 0.5102 & 1.5067 & 4.5686 & $\times 10^{+19}$ \\ 
$\pi/2$ & 3.606 & 9.977 & 0.757 & 4.525 & 8.257 & 10.094 & $\times 10^{+4}$\\ \hline
\end{tabular}
\end{center}
The conditional numbers for $\alpha=0$ increase exponentially with the increase of $k$. 
In fact, for values $k>250$ the determinant of $D$ has an infinitesimal magnitude of $e^{-196}$ and 
the double-precision float arithmetic is not enough for consistent calculations. 
On the other hand, the conditional number for $\alpha=\pi/2$ is rather stable and allows 
computations for $k$ well above $1000$.
In conclusion, the simulation results are in good support of the choice $\alpha = Ed(X,Y)$, $X,Y\sim Unif(U)$. 

The above described procedure can be used to improve the conditioning of a covariance operator. 
A way to achieve this is instead of definition (\ref{eq:rankmetric_covfield}), to use the following modified operator 
\begin{equation}\label{eq:rankmetric_covfield_alpha}
G\Sigma(p) := \int_{\mathcal{U}(p)} G(p)(Exp_p^{-1}q)(Exp_p^{-1}q)^T [1 - \frac{\alpha(p)}{d(p,q)}]^2 f(p) dV(p), 
\end{equation}
assuming that all $\alpha(p) := Ed^2(p,X)$, $X\sim Unif(\mathcal{U}(p))$ are well defined.
Theoretical and experimental analysis of (\ref{eq:rankmetric_covfield_alpha}), however, is out of the scope of this paper.

\begin{flushright}
\begin{minipage}[t]{0.5\linewidth}
		Department of Statistics, \\ 
		Florida State University, \\
		117 N. Woodward Ave., \\
        Tallahassee, Florida, 32306-4330 \\
		{\tt balov@stat.fsu.edu}
\end{minipage}
\end{flushright}

\end{document}